\documentclass{amsart} 

\usepackage{lineno, hyperref}

\usepackage{amsbsy,amsmath, amsfonts, latexsym, amsthm, amsxtra, amssymb, amssymb}
\usepackage[shortlabels]{enumitem}
\usepackage{soul}
\SetEnumerateShortLabel{a}{\textup{(\alph*)}}
\SetEnumerateShortLabel{A}{\textup{(\Alph*)}}
\SetEnumerateShortLabel{1}{\textup{(\arabic*)}}
\SetEnumerateShortLabel{i}{\textup{(\roman*)}}
\SetEnumerateShortLabel{(I)}{\textup{(\Roman*)}}

\newtheorem{thm}{Theorem}[section]

\newtheorem{lem}[thm]{Lemma}
\newtheorem{prop}[thm]{Proposition}

\theoremstyle{definition}

\newtheorem{df}[thm]{Definition}
\newtheorem{ex}[thm]{Example}

\newtheorem{ques}[thm]{Question}

\theoremstyle{remark}
\newtheorem{rem}[thm]{Remark}

\def\depth{\operatorname{depth}} 
\def\sdepth{\operatorname{sdepth}} 
\def\supp{\operatorname{supp}}

\def\NN{{\mathbb N}}
\def\ZZ{{\mathbb Z}}

\def\calD{{\mathcal D}}
\def\calF{ {\mathcal F}}
\def\calP{{\mathcal P}}

\def\frakm{{\mathfrak m}}

\def\ceil#1{\left\lceil #1 \right\rceil}
\def\floor#1{\left\lfloor #1 \right\rfloor}

\def\lcm{\operatorname{lcm}}

\def\Index#1{\emph{#1}}
\def\Set#1{\left\{ #1\right \}}

\title[Stanley depth of complete intersection monomial ideals]{Stanley depth of complete intersection monomial ideals and upper-discrete partitions}
\author[YiHuang Shen]{YiHuang Shen}
\address{Department of Mathematics, Purdue University, West Lafayette, IN 47907, U.S.A.} 
\email{yshen@math.purdue.edu}


\begin{document}

\begin{abstract}
    Let $I$ be an $m$-generated complete intersection monomial ideal in $S=K[x_1,\dots,x_n]$. We show that the Stanley depth of $I$ is $n-\floor{\frac{m}{2}}$.
We also study the upper-discrete structure for monomial ideals and prove that if $I$ is a squarefree monomial ideal minimally generated by $3$ elements, then the Stanley depth of $I$ is $n-1$.
\end{abstract}

\maketitle

\section{Introduction}
Let $\NN=\Set{0,1,2,\dots}$ be the set of non-negative integers. Let $K$ be a field and $S=K[x_1,\dots,x_n]$ be a polynomial ring over $K$. Suppose $M$ is a finitely generated $\ZZ^n$-graded $S$-module. If $u\in M$ is a homogeneous element and $Z$ is a subset of $\Set{x_1,\dots,x_n}$, then the $K$-subspace $uK[Z]$ of $M$ is called a \Index{Stanley space}. A \Index{Stanley decomposition of $M$} is a partition $\calD: M=\bigoplus_{i=1}^m u_i K[Z_i]$  in the category of $\ZZ^n$-graded $K$-vector spaces. The \Index{Stanley depth of $\calD$} is $\sdepth(\calD)=\min\Set{|Z_i|: 1\le i \le m}$ and the \Index{Stanley depth of $M$} is 
\[
\sdepth (M)=\max\Set{\sdepth(\calD): \text{ $\calD$ is  a Stanley decomposition of $M$}}.
\] 

The interest in finding Stanley decompositions and Stanley depths can be traced back to the pioneering paper of Stanley \cite{MR666158}.  There it was conjectured that $\depth(M) \le \sdepth(M)$. In \cite{MR2267659} it was shown that if $M$ allows a prime filtration $\calF$ with $\supp(\calF)=\min(M)$, then this conjecture holds. And if $I\subset S$ is a Gorenstein monomial ideal with $\dim(S) \le 5$, then \cite{MR2366164} showed that this conjecture is also true for $M=S/I$. However, in spite of the many supporting facts, the conjecture still remains open. One of the main obstacles for verifying the Stanley's conjecture lies in the difficulty of computing Stanley depths. Even with the method of Herzog, Vladoiu and Zheng which we will discuss immediately, it is still practically very difficult to find the Stanley depth for modules from general monomial ideals. The strongest result known to us that is pertinent to our work concerns the homogeneous maximal ideal $(x_1,\dots,x_n)\subset S$, which will be mentioned in Theorem \ref{bhkty} below.

In this paper, we will focus on the case where $M=I$ is a monomial ideal in $S$.
Let $G(I)=\Set{v_1,\dots,v_m}$ be the set of minimal monomial generators of $I$, and for $c=(c(1),\dots,c(n))\in \NN^n$, denote $x^c=\prod_i x_i^{c(i)}$. For a fixed $g\in \NN^n$ such that $\lcm(v_1,\dots,v_m)$ divides $x^g$, Herzog, Vladoiu and Zheng introduced in \cite{arxiv.0712.2308}  the associated poset $P_I^g=\Set{c\in \NN^n:  c\le g \text{ and $v_i |x^c$ for some $i$}}$ for $I$. Here $\le$ is the natural partial order in $\ZZ^n$ by componentwise comparison. For $a,b\in P_I^g$, define the interval $[a,b]$ to be $\Set{c\in P_I^g: a\le c\le b}$. Corresponding to each (disjoint) partition $\calP: P_I^g=\bigcup_{i=1}^r [c_i, d_i]$, there is a Stanley decomposition $\calD(\calP)$ of $I$.   They showed in \cite[Corollary 2.5]{arxiv.0712.2308} that there is a partition $\calP$ such that $\sdepth(I)=\sdepth(\calD(\calP))$.

Recently, Cimpoea\c{s} studied Stanley decomposition of complete intersection ideals. He proved in \cite[Theorem 2.1]{arxiv.0805.2306} that the Stanley depth of a complete intersection monomial ideal is equal to the Stanley depth of its radical.  Therefore, the focus of research is directed to squarefree monomial ideals.  Recall that a Stanley space $uK[Z]$ is called squarefree, if $u$ is squarefree and $\supp(u)\subset Z$. 
If $I$ is a squarefree monomial ideal, we can take $g=(1,\dots,1)$ and  write $P_I^g$ simply as $P_I$. Recall that a vector $d\in \ZZ^n$ is \Index{squarefree} if $d(i)=0$ or $1$, for all $1\le i \le n$. If $d\in \NN^n$ is squarefree, write $Z_d=\Set{x_j: d(j)=1}$. Then for any partition  $\calP: P_I=\bigcup_i [c_i,d_i]$, $\calD(\calP): I=\bigoplus_i x^{c_i}K[Z_{d_i}]$ is the associated Stanley decomposition of $I$ introduced in \cite{arxiv.0712.2308}. Meanwhile $\sdepth(\calD(\calP))=\min\Set{|d_i|: 1\le i \le r}$. Here $|d_i|$ is the sum of components in $d_i$.  The Stanley decomposition $\calD(\calP)$ is clearly squarefree. This observation shows in particular that
\[
\sdepth(I)=\max\Set{\sdepth(\calD): \text{$\calD$ is a squarefree Stanley decomposition of $I$}}.
\]

This paper proceeds as follows. We compute in Theorem \ref{ci} the Stanley depth of complete intersection monomial ideals. It turns out that the Stanley depth depends only on the dimension of the polynomial ring and the minimal number of generators. The third section studies the upper-discrete partition of squarefree monomial ideals. And in the last section, we prove that the Stanley depth of a squarefree monomial ideal minimally generated by $3$ elements is $n-1$. For $4$-generated squarefree monomial ideals, the lower bound of Stanley depth is $n-2$.

\section{Stanley depth of complete intersection monomial ideals}
The Stanley depth of the monomial maximal ideal is known.
\begin{thm}
    [{\cite[Theorem 2.2]{BHKTY}}] 
    \label{bhkty}
    Let $\frakm=(x_1,\dots,x_n)$ be the maximal ideal in $S=K[x_1,\dots,x_n]$, then  $\sdepth(\frakm)=\ceil{\frac{n}{2}}$.
\end{thm}

Herzog, Vladoiu and Zheng computed the Stanley depth of 3-generated  complete intersection monomial ideals.

\begin{prop}
    [{\cite[Proposition 3.8]{arxiv.0712.2308}}]
    Let $I\subset S$ be a  complete intersection monomial ideal minimally generated by $3$ elements. Then $\sdepth(I)=n-1$.
\end{prop}

We want to generalize the above two results and answer Conjecture 2.5 in \cite{arxiv.0805.2306}.  For simplicity of notation, we identify any squarefree vector $c\in \ZZ^n$ with $\Set{i\mid c(i)=1}$. 

\begin{lem}
    \label{LEM}
    Let $v_1,\dots,v_m$ be  squarefree monomials in $K[x_1,\dots,x_{n-1}]$. If $I=(v_1,\dots,v_{m-1},v_m x_{n})$ and $I'=(v_1,\dots,v_{m-1},v_m x_{n}x_{n+1})$ are ideals in $S=K[x_1,\dots,x_n]$ and $S'=S[x_{n+1}]$ respectively, then $\sdepth(I')=\sdepth(I)+1$.
\end{lem}

\begin{proof}
    By assumption, there is a partition $\calP: P_I=\bigcup_i [c_i,d_i]$ for $I$ such that $\sdepth(\calD(\calP))=\sdepth(I)$.  By   \cite[Corollary 2.3]{arxiv.0805.2306}, $\sdepth(I')\le \sdepth(I)+1$. Now it suffices to construct a partition $\calP'$ for $P_{I'}$ with $\sdepth(\calD(\calP'))=\sdepth(I)+1$.
    
    For each interval $B=[c,d]$ in $\calP$, we define the corresponding interval $B'$:
    \begin{enumerate}[1]
        \item If $n\in c$, which by our identification means $c(n)=1$,  let $B^1=[c\cup \Set{n+1},d\cup \Set{n+1}]$.
        \item If $n\not \in c$, let $B^2=[c,d\cup \Set{n+1}]$. Furthermore, if $n\not\in d$, let $B^3=[c\cup \Set{n},d\cup \Set{n}]$.
    \end{enumerate}
    Let $B'$ be the union of those $B^k$'s defined. Hence $B'=B^1$, $B'=B^2$ or $B'=B^2 \cup B^3$. $B'$ is a subset of $P_{I'}$. We claim that $\calP': P_{I'}=\bigcup_{i=1}^r B_i'$ is a partition for $P_{I'}$ with $\sdepth(\calD(\calP'))=\sdepth(I)+1$.

    First, we prove that the intervals $B_i'$ cover $\calP'$. Let $u$ be a proper subset of $\Set{1,\dots, n+1}$ in $P_{I'}$.  Depending on whether $n+1 \in u$, we have two cases.

    \begin{enumerate}[1]
        \item If $n+1 \in u$, let $u'=u\setminus \Set{n+1}$. We have $u'\in P_I$, hence there is an interval  $B=[c,d]$ in $\calP$ such that $u'\in B$. If $n\in c$, then $u\in B^1$. Otherwise, $n\not\in c$, and  $u\in B^2$.
        \item If $n+1 \not \in u$, then $x^u$ is divisible by some $v_i\ne v_m$. Consequently, $u\in P_I$ and there is an interval  $B=[c,d]$ in $\calP$ with $u\in B$. 
                    \begin{enumerate}[a]
                        \item If $n\not\in c$, then $u\in B^2$.
                        \item If $n\in c$, then $n\in u$ as well. Let $u'=u\setminus\Set{n}$ and again we have $u'\in P_I$. There is an interval  $\tilde{B}=[\tilde{c},\tilde{d}]$ in $\calP$ with $u'\in \tilde{B}$. Since $n\not\in u'$, $n\not\in \tilde{c}$.  Now depending on whether $n\in \tilde{d}$ or $n\not \in \tilde{d}$, $u\in \tilde{B}^2$ or $u\in \tilde{B}^3$.
            \end{enumerate}
    \end{enumerate}

    Now we show that the intervals in $\calP'$ are pairwise disjoint. Suppose  $B_1=[c_1,d_1]$ and $B_2=[c_2,d_2]$ are intervals in $\calP$. We prove by contradiction that $B_1^i$ and $B_2^j$ are  disjoint for $1\le i\ne j \le 3$.

    Suppose that $u\in B_1^1 \cap B_2^2$, then $n+1\in u$. Let $u'=u\setminus\Set{n+1}$. Then $u'\in B_1 \cap B_2$, hence $B_1=B_2$. But $n\in c_1$,  $n\not\in c_2$ and $c_1=c_2$. This is a contradiction.

    Suppose that $u\in B_1^1 \cap B_2^3$, then $n+1 \in u$. But $n+1\not\in d$, and $x^{d\cup\Set{n}}$ is divisible by $x^u$. As a result, $n+1 \not \in u$. This is a contradiction.

    Suppose that $u\in B_1^2 \cap B_2^3$, then $n\in u$ and $n+1 \not \in u$. Let $u'=u\setminus \Set{n}$, then $u'\in B_2$. Since $n\not\in c_1$ and $x^u$ is divisible by $x^{c_1}$,  $x^{u'}$ is also divisible by $x^{c_1}$. Meanwhile, since $n+1 \not\in u$,  $n+1\not\in u'$. Since $x^{d_1\cup\Set{n+1}}$ is divisible by $x^u$, $x^{d_1}$ is divisible by $x^{u'}$. Thus $u'\in B_1$ as well. Hence $B_1=B_2$. Now since $u\in B_2^3$, $n\in u$ and $n\not \in d_2$. Since $d_1=d_2$ and $u\in B_1^2$, $n\not\in u$. This is a contradiction.

    Now let $i=j$. If $u\in B_1^1\cap B_2^1$ or $B_1^2 \cap B_2^2$, let $u'=u\setminus\Set{n+1}$, then $u\in B_1\cap B_2$ and $B_1=B_2$. Likewise, if $u\in B_1^3\cap B_2^3$, let $u'=u\setminus\Set{n}$. Notice that $n+1\not\in u'$ and $u'\in B_1\cap B_2$. Hence $B_1=B_2$.

\end{proof}

\begin{thm}
    \label{ci}
    Let $I\subset S=K[x_1,\dots,x_n]$ be a  complete intersection monomial ideal minimally generated by $m$ elements. Then $\sdepth(I)=n-\floor{\frac{m}{2}}$.
\end{thm}

\begin{proof}
    Following \cite[Theorem 2.1]{arxiv.0805.2306} and \cite[Lemma 3.6]{arxiv.0712.2308}, one can assume that $I$ is squarefree and every ring variable shows up in exactly one monomial generator of $I$. We fix $m$ and prove the theorem by induction on $n\ge m$. 
   
    The base case is when $n=m$ and  hence $I=(x_1,\dots,x_m)$ is the maximal ideal.  The validity now follows from Theorem \ref{bhkty}. Notice that $\ceil{\frac{m}{2}}=m-\floor{\frac{m}{2}}$.

    Now let $n\ge m$ and assume that the theorem holds for $n$. We want to prove that it also holds for $n+1$. Without loss of generality, we consider a squarefree complete intersection monomial ideal $I'$ in $ S'=S[x_{n+1}]$, minimally generated by monomials $v_1,\dots,v_{m-1}, v_m x_{n+1}$ and assume that $x_n$ divides $v_m$. Then the ideal $I=(v_1,\dots,v_m)$ in $S$ is also a squarefree complete intersection monomial ideal. Therefore, by the induction hypothesis, $\sdepth(I)=n-\floor{\frac{m}{2}}$. Now by Lemma \ref{LEM}, $\sdepth(I')=\sdepth(I)+1$ and this completes the proof.
\end{proof}

\section{Upper-discrete partitions}
In this section, we introduce the upper-discrete partitions. It will be the main tool in the next section to study the 3-generated squarefree monomial ideals.

\begin{df}
    Let $P$ be the associated poset of monomials in $S$. $P$ is called \Index{upper-discrete of degree $k$}, if there is a partition $\calP: P=\bigcup_i [c_i,d_i]$, such that $|d_i|\ge k$ for all $i$, and $c_i=d_i$ when $|d_i|>k$. And this partition is called an \Index{upper-discrete partition of degree $k$}. 
\end{df}

\begin{ex}
    We use the notations in figure 2 of \cite{arxiv.0712.2308} and consider the ideal $I=(x_1x_2,x_2x_3,x_1x_3)\subset K[x_1,x_2,x_3]$. It is readily seen that $P_I = [12,12]\cup [23,23]\cup [13,13]\cup [123,123]$  gives an upper-discrete partition of degree $2$. However, a shorter one $P_I=[12,123]\cup [23,23]\cup [13,13]$ does not.
\end{ex}

\begin{prop}
If $I$ is a squarefree monomial ideal in $S$,  then the poset $P_I$ is upper-discrete of degree $k$ for $k\le \sdepth(I)$.
\end{prop}

\begin{proof}
    Let $\calP:P_I=\bigcup_i [c_i,d_i]$ be a Stanley decomposition with $\sdepth(\calD(\calP))=\sdepth(I)$. Hence $|d_i|\ge \sdepth(I)\ge k$. Now it suffices to show that each interval $[c_i,d_i]$ allows an upper-discrete partition of degree $k$. This is equivalent to say that the interval $[\emptyset, d_i\setminus{c_i}]$ admits an upper-discrete partition of degree $k-|c_i|$, where $[\emptyset, d_i\setminus{c_i}]$ is an interval in the poset $P_{S}$ for the unit ideal $S$. Since $[\emptyset, d_i\setminus{c_i}]$ is isomorphic to the poset $P_{S'}$ where $S'=K[x_1,\dots,x_{|d_i|-|c_i|}]$, it is enough to show that the poset $P_{S}$ is upper-discrete of degree $k$ for $0\le k \le n$.

    We prove by induction on $n$. The base cases when $n=0$ or $n=1$ are trivial. Now let $n\ge 2$ and suppose the claim holds for $n-1$. The cases when $k=0$ or $k=n$ are clear. Hence  we may assume that $1\le k \le n-1$ and let $S'=K[x_1,\dots,x_{n-1}]$. Then  $P_{S'}$ have two upper-discrete partitions $\calP^1: P_{S'}=\bigcup_i [c_i^1,d_i^1]$ and $\calP^2: P_{S'}=\bigcup_i [c_i^2,d_i^2]$ of degrees $k$ and $k-1$ respectively. Clearly 
    \[\calP:P_{S}=(\bigcup_i[c_i^1,d_i^1]) \cup (\bigcup_i[c_i^2\cup\Set{n},d_i^2\cup \Set{n}])
    \]
    is an upper-discrete partition of degree $k$. And this completes the proof.
\end{proof}

\begin{rem}
    \label{REM}
    Let $I\subset S$ be a squarefree complete intersection monomial ideal with minimal monomial generating set $G(I)=\Set{v_1,\dots,v_m}$. We further  assume that $x_n$ divides $v_m$. Let  $I'=(v_1,\dots,v_m x_{n+1})\subset S'=S[x_{n+1}]$. If $P_I$ has an upper-discrete partition $\calP$ of degree $k$, then the proof of Lemma \ref{LEM} can be modified as follows to give an upper-discrete partition of $P_{I'}$ of degree $k+1$.
    
    Let $B=[c,d]$ be an interval in $\calP$. We construct the interval $B'$ in the following way:
       \begin{enumerate}[1]
        \item If $n\in c$, let $B^1=[c\cup \Set{n+1},d\cup \Set{n+1}]$.
        \item If $n\not \in c$,
            \begin{enumerate}[a]
                \item if $|c|\le k$, let $B^2=[c,d\cup \Set{n+1}]$.  Furthermore, if $n\not\in d$, let $B^3=[c\cup \Set{n},d\cup \Set{n}]$.
                \item if $|c|> k$, hence $c=d$, then let 
                    \begin{itemize}
                        \item  $B^4=B$, 
                        \item  $B^5=[c\cup\Set{n},c\cup\Set{n}]$, 
                        \item  $B^6=[c\cup\Set{n+1},c\cup\Set{n+1}]$.
                    \end{itemize}
                             \end{enumerate}
    \end{enumerate}
    Let $B'$ be the union of those $B^k$ defined. Hence either $B'=B^1$, $B'=B^2$, $B'=B^2 \cup B^3$, or $B'=B^4\cup B^5 \cup B^6$. The rest of the proof is essentially the same.
\end{rem}

\section{Squarefree monomial ideals}

If $I$ is not a complete intersection, the formula in Theorem \ref{ci} will fail in general. For instance, let $I=(x_1 x_2 x_3,x_1x_2x_4,x_1x_3x_4,x_2x_3x_4)$ in $S=K[x_1,\dots,x_4]$. Then $\sdepth(I)=3$ instead of $4-\floor{\frac{4}{2}}=2$. However, when $m=3$, the situation is different. 

\begin{thm}
    \label{3GEN}
    Let $I$ be a 3-generated squarefree monomial ideal in $S=K[x_1,\dots,x_n]$. Then $\sdepth(I)\ge n-1$. In particular, if $I$ is not principal, $\sdepth(I)=n-1$.
\end{thm}

\begin{proof}
    Let $I$ be  generated by monomials $v_1$, $v_2$ and $v_3$. For any ring variable $x_j$, we say $x_j$ is of type $i$, if there are exactly $i$ of the three generators involve the variable $x_j$.

    If $x_n$ is of type $0$, then for the ideal $I'=(v_1,v_2,v_3)$ in $K[x_1,\dots,x_{n-1}]$, we have $\sdepth(I')=\sdepth(I)-1$ by \cite[Lemma 3.6]{arxiv.0712.2308}.

    In a like manner, if $x_n$ is of type $3$, then for the ideal $I'=(v_1/x_n,v_2/x_n,v_3/x_n)$ in $S$, it is readily seen that $I'$ is naturally isomorphic to $I$ in the category of $\ZZ^n$-graded $K$-vector spaces up to degree shifting. Thus, $\sdepth(I')=\sdepth(I)$. But then, $x_n$ is of type $0$ for $I'$.

    Hence it suffices to prove the result for the case when all ring variables are of type either $1$ or $2$. We call variable $x_j$ to be of type $1$-$(i)$, if $x_j$ is of type $1$ and $v_i$ involves $x_j$.  By Lemma \ref{LEM}, we may assume that for every $i$, $1\le i\le 3$, there is at most one ring variable to be of type $1$-$(i)$. 

After these reductions, it is easily seen that the proof is done once we can show the following.
\begin{enumerate}[(I)]
\item   Fix $n\ge 0$ and  let $I$ be any  ideal in $S=K[x_1,\dots,x_n]$ generated by squarefree monomials $v_1$, $v_2$ and $v_3$, such that all ring variables are of type $2$ for $I$. We prove that $\sdepth(I)\ge n-1$.
 \item For any fixed $I$ in (I), we also consider ideals $I_1=(v_1x_{n+1},v_2,v_3)$ in $S_1=S[x_{n+1}]$, $I_2=(v_1x_{n+1},v_2x_{n+2},v_3)$ in $S_2=S_1[x_{n+2}]$, and $I_3=(v_1x_{n+1},v_2 x_{n+2}, v_3 x_{n+3})$ in $S_3=S_2[x_{n+3}]$. We prove that $\sdepth(I_i)\ge n-1+i$ for $1\le i \le 3$.
\end{enumerate}
   
The proof is then carried out in 4 steps.

\paragraph{Step 0}

   To begin with, we investigate the ideal $I$ in case (I) and assume that all ring variables are of type $2$. We prove by induction on $n\in \NN$ that $\sdepth(I)\ge n-1$. 
   
   The base cases when $n\le 1$ are easy to verify. Now we assume that the formula holds for a fixed $n\ge 1$ and consider the ideal $I'=(v_1x_{n+1},v_2x_{n+1},v_3)$ in $S'=S[x_{n+1}]$. Here $v_1$, $v_2$ and $v_3$ are squarefree monomials in $S=K[x_1,\dots,x_n]$,  and all ring variables of $S$ are of type $2$ for $I=(v_1,v_2,v_3)$ in $S$. We want to show that $\sdepth(I')\ge n$.
   
   By induction hypothesis, $\sdepth(I)\ge n-1$. Thus we can find $\calP: P_I=\bigcup_i [c_i,d_i]$,  an upper-discrete partition of degree $n-1$. For each interval $B=[c,d]$ in $\calP$, define $B'$ as follows.
\begin{enumerate}[1]
    \item Suppose $|d|= n-1$. If $v_3$ divides $x^c$, let $B^1=[c, d\cup\Set{n+1}]$. Otherwise, $v_3 \nmid x^c$, and let $B^2=[c\cup\Set{n+1},d\cup\Set{n+1}]$.
    \item If $|d|=n$, then $c=d=\Set{1,\dots,n}$. Let $B^3=[c,c]$ and $B^4=[c\cup\Set{n+1},c\cup\Set{n+1}]$.
\end{enumerate}
Let $B'$ be $B^1$, $B^2$ or $B^3\cup B^4$ correspondingly. $B'$ is a subset of $P_{I'}$. We claim that $\calP': P_{I'}=\bigcup B_i'$ is an upper-discrete partition of degree $n$.

We first show that the intervals $B_i'$ cover $P_{I'}$. Let $u\in P_{I'}$, if $u=\Set{1,\dots,n+1}$, then $u\in B^4$. Otherwise, we may assume that $|u|\le n$.
\begin{enumerate}[1]
    \item If $n+1\not\in u$, then $v_3$ divides $u$. For this reason, we have $u\in P_I$, and $u\in B=[c,d]$ in $\calP$. If $u=\Set{1,\dots,n}$, then $u\in B^3$. Otherwise, $|u|\le n-1$.  We claim that $v_3 \mid x^c$, hence $u\in B^1$. If $v_3 \nmid x^c$, we have  $v_1$ or $v_2$ dividing $x^c$. As a result, $x^d$ is divisible by $v_1$ or $v_2$. But $x^d$ is also divisible by $v_3$, thus $\supp(v_1 v_3)=\supp(v_2 v_3)=\Set{1,\dots,n}\subset d$. At the same time, $|d|=n-1$ and this is a contradiction.

    \item If $n+1\in u$, let $u'=u\setminus\Set{n+1}$ and we  have $u'\in P_I$. Hence there is an interval $B=[c,d]$ in $\calP$ with $u'\in B$. Then depending on whether $v_3 \mid x^c$ or $v_3 \nmid x^c$, $u\in B^1$ or $u\in B^2$.
\end{enumerate}

Now we show the intervals in $\calP'$ are pairwise disjoint. It is straightforward to check that $\Set{1,\dots,n}$ and $\Set{1,\dots,n+1}$ are only in intervals $B^3$ and $B^4$ respectively. 

Now suppose $u\in B_1^1 \cap B_2^2 \ne \emptyset$ for $B_1=[c_1,d_1]$ and $B_2=[c_2,d_2]$ in $\calP$. According to the construction, $v_3 \mid x^{c_1}$ and $c_1\le u \le d_2\cup\Set{n+1}$. Hence $v_3 \mid x^{d_2\cup\Set{n+1}}$. At the same time, $v_3\nmid x^{c_2}$, hence $v_1\mid x^{c_2}$ or $v_2\mid x^{c_2}$.  Thus $x^{d_2}$ is divisible by either $v_1$ or $v_2$, and $d_2=\Set{1,\dots,n}$. This is against the assumption that $|d_2|=n-1$.

Likewise, if $u\in B_1^1\cap B_2^1$ or $B_1^2 \cap B_2^2$, then $u\setminus\Set{n+1}\in B_1 \cap B_2$. Hence $B_1=B_2$. This completes the proof for the claim.

\paragraph{Step 1}
Let $I$ be the ideal in case (I). Then in Step 0, we showed that $\sdepth(I)\ge n-1$. As a result,  we have  an upper-discrete partition  $\calP:P_I=\bigcup_i [c_i,d_i]$  of degree $n-1$. Now we construct an upper-discrete partition of degree $n$ for $P_{I_1}$, where $I_1$ is constructed in case (II).

    For each $B=[c,d]$ in $\calP$,  we define $B'$ as follows. 
    \begin{enumerate}[1]
        \item Suppose $|d|=n-1$.  If $v_1$ divides $x^c$, then let $B^1=[c\cup \Set{n+1},d\cup\Set{n+1}]$. Otherwise, $v_1 \nmid x^c$, and let $B^2=[c,d\cup\Set{n+1}]$. 
        \item Suppose $|d|=n$, then $c=d=\Set{1,\dots,n}$. Let $B^3=[c,c]$ and $B^4= [c\cup\Set{n+1},c\cup\Set{n+1}]$.
    \end{enumerate}
    Define $B'=B^1$, $B'=B^2$ or $B'=B^3 \cup B^4$ correspondingly. $B'$ is a subset of $P_{I_1}$. We claim that $\calP_1: P_{I_1}=\bigcup_i B_i'$ is a partition that satisfies the requirement. 

     We first show that intervals $B_i'$  cover $P_{I_1}$.  Let $u\in P_{I_1}$. If $|u|=n+1$, then $u=\Set{1,\dots,n+1}$, and $u\in B^4$.  Otherwise, we may assume that $|u|\le n$.
     \begin{enumerate}[1]
        \item If $v_1 x_{n+1}$ divides $x^u$, then $n+1 \in u$. Let $u'=u\setminus\Set{n+1}$, then $v_1$ divides $x^{u'}$ and $u'\in P_I$. Thus there is an interval $B=[c,d]$ in $\calP$ such that $u'\in B$. Since $|u'|\le n-1$, $|d|=n-1$. We claim that $v_1\mid x^c$, hence $u\in B^1$. Otherwise, $v_2$ or $v_3$ divides $x^c$. Say it is $v_2$, then $v_2$ also divides $x^d$. On the other hand, $v_1$ divides $x^{u'}$, hence $x^d$ is also divisible by $v_1$. Thus $\supp(v_1 v_2)=\Set{1,\dots,n}\subset d$ and  $|d|\ge n$. However $|d|=n-1$ and this is a contradiction.

        \item If $v_1$ does not divide $x^u$, then neither does $v_1 x_{n+1}$. Therefore, $v_2$ or $v_3$ divides $x^u$. Let $u'=u\setminus\Set{n+1}$. Then $v_2$ or $v_3$  divides $x^{u'}$, and we have $u'\in P_I$.  Let $u'\in B=[c,d]$, an interval in $\calP$. Since  $v_1 \nmid x^{u'}$, we have $v_1 \nmid x^c$ and $u\in B^2$.
        \item If $v_1$ divides $x^u$, but $v_1 x_{n+1} $ does not, then $n+1 \not\in u$. Since $u\in P_{I_1}$,  $x^u$ is divisible by $v_2$ or $v_3$. Since $\supp(v_1v_2)=\supp(v_1v_3)=\Set{1,\dots,n}$, this would force $u=\Set{1,\dots,n}$, and $u\in B^3$.
    \end{enumerate}

    Now we show that $\calP_1: P_{I_1}=\bigcup_i B_i'$ is a disjoint union. Since $\calP$ is an upper-discrete partition, if $u_1=\Set{1,\dots,n+1}$, $B^4$ is the only interval containing $u_1$. 

    Consider $u_2=\Set{1,\dots,n}$ and suppose that $u_2\in B^1$ for some $B=[c,d]$ in $\calP$. Then $n+1\in u_2$ and this is impossible. On the other hand, suppose  $u_2\in B^2$ for some $B=[c,d]$ in $\calP$. Then $c\le u_2 \le d\cup\Set{n+1}$. Since $n+1\not\in u_2$, we have $c\le u_2\le d$ and $u_2\in B$. Hence $|d|\ge |u_2|=n$. On the other hand, by our assumption on $B^2$, $|d|=n-1$ and this is a contradiction. 

    Let $B_1=[c_1,d_1]$ and $B_2=[c_2,d_2]$ be intervals in $\calP$. If $u\in B_1^1 \cap B_2^2 \ne \emptyset$, then $u\setminus\Set{n+1}\in B_1 \cap B_2$. Hence $B_1=B_2$. Meanwhile, $v_1 \mid x^{c_1}$, while $v_1 \nmid x^{c_2}$. This is a contradiction.

    Similarly, if $u\in B_1^1 \cap B_2^1$ or $B_1^2 \cap B_2^2$, then $u\setminus\Set{n+1}\in B_1 \cap B_2$. Thus $B_1=B_2$.

    \paragraph{Step 2}

    Using partition $\calP_1$ in Step 1, we construct an upper-discrete partition $\calP_2$ for $P_{I_2}$ with degree $n+1$.  For each $B=[c,d]$ in $\calP_1$,  we define $B'$ as follows. 
    \begin{enumerate}[1]
        \item Suppose $|d|=n$.  If $v_2$ divides $x^c$, then let $B^1=[c\cup \Set{n+2},d\cup\Set{n+2}]$. Otherwise, $v_2 \nmid x^c$, and let $B^2=[c,d\cup\Set{n+2}]$. 
        \item Suppose $|d|=n+1$, then $c=d=\Set{1,\dots,n+1}$. Let $B^3=[\Set{1,\dots,n},c]$ and $B^4= [c\cup\Set{n+2},c\cup\Set{n+2}]$.
    \end{enumerate}
    Define $B'=B^1$, $B'=B^2$ or $B'=B^3 \cup B^4$ correspondingly. $B'$ is a subset of $P_{I_2}$. We claim $\calP_2: P_{I_2}=\bigcup_i B_i'$ is a partition that satisfies the requirement. 

  We first show that intervals $B_i'$  cover $P_{I_2}$.  Let $u\in P_{I_2}$. If $|u|=n+2$, then $u=\Set{1,\dots,n+2}$, and $u\in B^4$.  Otherwise, we may assume that $|u|\le n+1$.
     \begin{enumerate}[1]
         \item If $v_2 x_{n+2}$ divides $x^u$, then $n+2 \in u$. Let $u'=u\setminus\Set{n+2}$, then $v_2$ divides $x^{u'}$ and $u'\in P_{I_1}$. Thus there is an interval $B=[c,d]$ in $\calP_1$ that $u'\in B$. Since $|u'|\le n$, $|d|=n$. We claim that $v_2\mid x^c$, hence $u\in B^1$. Otherwise, $v_1 x_{n+1}$ or $v_3$ divides $x^c$. 
            
             \begin{enumerate}
                 \item   If $v_1 x_{n+1}$ divides $x^c$, then $x^d$ is divisible by both $v_1 x_{n+1}$ and $v_2$. Thus $\supp(v_1 x_{n+1} v_2)=\Set{1,\dots,n+1}\subset d$ and  $|d|\ge n+1$. Nevertheless, $|d|=n$ and this is a contradiction.

                 \item   If $v_3$ divides $x^c$, then $x^d$ is divisible by both $v_2$ and $v_3$. Thus $\supp(v_2 v_3)=\Set{1,\dots,n}\subset d$. Since $|d|=n$, $d=\Set{1,\dots,n}$. Thus by the construction of $\calP_1$, $c=d$. We still have $v_2 \mid x^c$.
             \end{enumerate}
           
         \item If $v_2$ does not divide $x^u$, then neither does $v_2 x_{n+2}$. Hence $v_1 x_{n+1}$ or $v_3$ divides $x^u$. Let $u'=u\setminus\Set{n+2}$. Then $v_1 x_{n+1}$ or $v_3$  divides $x^{u'}$, and we have $u'\in P_{I_1}$.  Let $u'\in B=[c,d]$, an interval in $\calP_1$. Since  $v_2 \nmid x^{u'}$, we have $v_2 \nmid x^c$ and $u\in B^2$.
         \item If $v_2$ divides $x^u$, but $v_2 x_{n+2} $ does not, then $n+2 \not\in u$. Thus $x^u$ is divisible by $v_1 x_{n+1}$ or $v_3$. 
             \begin{enumerate}
                 \item  If $x^u$ is divisible by $v_1 x_{n+1}$, since $\supp(v_1 x_{n+1} v_2)=\Set{1,\dots,n+1}$, this would force $u=\Set{1,\dots,n+1}$, and $u\in B^3$.
                 \item Otherwise, $x^u$ is divisible by $v_3$. At this moment, $\supp(v_2 v_3)=\Set{1,\dots,n}\subset u$. Since $|u|\le n+1$ and $n+1, n+2 \not\in u$, this forces $u=\Set{1,\dots,n}$, and $u\in B^3$.
             \end{enumerate}
                \end{enumerate}

    Now we need to show that $\calP_2: P_{I_2}=\bigcup_i B_i'$ is a disjoint union. The proof is similar to that in Step 1. However, one still need to consider $u_3=\Set{1,\dots,n}$.

    If $u_3 \in B^1$ for some $B=[c,d]\in \calP_1$, then $n+2 \in u_3$. This is impossible. If $u_3\in B^2$, then $c\le u_3 \le d\cup \Set{n+2}$. Since $n+2 \not\in u_3$, this implies that $c\le u_3 \le d$, i.e., $u_3 \in B$. By our construction of $\calP_1$, $c=d$. Thus $v_2$ divides $x^c$, and instead of $B^2$, we should construct $B^1$. This is a contradiction.

\paragraph{Step 3}

Using partition $\calP_2$ in Step 2, we construct an upper-discrete partition $\calP_3$ for $P_{I_3}$ with degree $n+2$.  For each $B=[c,d]$ in $\calP_2$,  we define $B'$ as follows. 
    \begin{enumerate}[1]
        \item  Suppose $\Set{1,\dots,n}\not\subset d$, then $|d|=n+1$. If $v_3$ divides $x^c$, then let $B^1=[c\cup \Set{n+3},d\cup\Set{n+3}]$. Otherwise, $v_3 \nmid x^c$, and let $B^2=[c,d\cup\Set{n+3}]$. 

        \item Suppose $\Set{1,\dots,n}\subset d$.  Then according to the construction of $\calP_2$, $B$ is one of the following intervals: 
            \begin{itemize}
                \item $[\Set{1,\dots,n},\Set{1,\dots,n+1}]$, 
                \item $[\Set{1,\dots,n,n+2},\Set{1,\dots,n,n+2}]$, 
                \item $[\Set{1,\dots,n+2},\Set{1,\dots,n+2}]$.
            \end{itemize}
            In particular, $\Set{1,\dots,n}\subset c$. Now define 
            \[
            B^3=[\Set{1,\dots,n+1},\Set{1,\dots,n+2}],
            \]
            \[
            B^4=[\Set{1,\dots,n,n+2},\Set{1,\dots,n,n+2,n+3}],
            \]
\[
B^5=[\Set{1,\dots,n,n+3},\Set{1,\dots,n+1,n+3}],
\]
            and 
            \[
            B^6=[\Set{1,\dots,n+3},\Set{1,\dots,n+3}.
            \]
    \end{enumerate}
    Define $B'=B^1$, $B'=B^2$ or $B'=B^3 \cup B^4\cup B^5 \cup B^6$ correspondingly. $B'$ is a subset of $P_{I_3}$.  We claim $\calP_3: P_{I_3}=\bigcup_i B_i'$ is a partition that satisfies the requirement. 

    We first show that intervals $B_i'$  cover $P_{I_3}$. Let $u\in P_{I_3}$. If $\Set{1,\dots,n}\subset u$, then $|u|\ge n+1$ and $u$ is in exactly one of the $B^i$ for $3\le i \le 6$. Otherwise, we have $|u|\le n+2$ and $x^u$ is divisible by exactly one of the monomial generators $v_i x_{n+i}$ for $I_3$.
     \begin{enumerate}[1]
         \item If $v_3 x_{n+3}$ divides $x^u$, then $n+3 \in u$. Let $u'=u\setminus\Set{n+3}$, then $v_3$ divides $x^{u'}$ and $u'\in P_{I_2}$. Thus there is an interval $B=[c,d]$ in $\calP_2$ that $u'\in B$. If $v_1 x_{n+1}$ or $v_2 x_{n+2}$ divides $x^c$, then $x^u$ is also divisible by it, which is impossible. Hence $v_3\mid x^c$. Since $\Set{1,\dots, n} \not\subset u$, by our construction of $\calP_2$, $\Set{1,\dots,n}\not\subset d$ and $u\in B^1$.

         \item If $v_3$ does not divide $x^u$, then neither does $v_3 x_{n+3}$. Hence $v_1 x_{n+1}$ or $v_2 x_{n+2}$ divides $x^u$. Let $u'=u\setminus\Set{n+3}$. Then $v_1 x_{n+1}$ or $v_2 x_{n+2}$  divides $x^{u'}$, and we have $u'\in P_{I_2}$.  Let $u'\in B=[c,d]$, an interval in $\calP_2$. Since  $v_3 \nmid x^{u'}$, we have $v_3 \nmid x^c$ and $u\in B^2$.
         \item If $v_3$ divides $x^u$, but $v_3 x_{n+3} $ does not, then $n+3 \not\in u$. Thus $x^u$ is divisible by $v_1 x_{n+1}$ or $v_2 x_{n+2}$. Hence $\supp(v_1 v_3)=\supp(v_2 v_3)=\Set{1,\dots,n}\subset u$, which is a contradiction.
             \end{enumerate}

             At this stage, we have to show that $\calP_2: P_{I_2}=\bigcup_i B_i'$ is a disjoint union. Let $u\in P_{I_3}$. If $\Set{1,\dots,n}\subset u$, then $|u|\ge n+1$ and $u\in B^3\cup B^4 \cup B^5 \cup B^6$. Suppose $u\in \tilde{B}^1$ or $\tilde{B}^2$ for some interval $\tilde{B}=[c,d]$ in $\calP_2$. Then $\Set{1,\dots,n}\subset d\cup\Set{n+3}$. This implies that $\Set{1,\dots,n}\subset d$, which is against the construction of $\tilde{B}^1$ or $\tilde{B}^2$. The rest of the proof is similar to that in Step 1. 
\end{proof}

As shown by the example at the beginning of this section, the Stanley depth of 4-generated squarefree monomial ideal is not necessarily $n-2$. Nevertheless, $n-2$ is the sharp lower bound.

\begin{prop}
    \label{PROP4}
    Let $I\subset S=K[x_1,\dots,x_n]$ be a squarefree monomial ideal generated by $4$ elements. Then $\sdepth(I)\ge n-2$.
\end{prop}

\begin{proof}
    We apply the technique in \cite[Proposition 3.4]{arxiv.0712.2308} and use their notations. Hence we prove by induction on $n$, with $n=1$ being trivial. Now consider $n\ge 2$ and assume that the claim holds for $n-1$. Suppose the minimal monomial generating set is $G(I)=\Set{x^{a_1},\dots,x^{a_4}}$. Without loss of generality, we may assume that $a_1 \vee \cdots \vee a_4=(1,\dots,1)$. Then there is a disjoint union $P_I=A_0 \cup A_1$, where $A_i=\Set{c\in P_I: c(n)=i}$ for $0\le i \le 1$.

    It is observed in \cite{arxiv.0712.2308} that $A_i=\Set{(c,i): c\in P_{I_i}^g}$ with $g=(1,\dots,1)\in \NN^{n-1}$, and $I_i$ is the monomial ideal in $K[x_1,\dots,x_{n-1}]$ such that 
    \[
    I \cap x_n^i K[x_1,\dots,x_{n-1}]=x_n^i I_i.
    \]
    $I_0$ and $I_1$ are still squarefree. Furthermore, $|G(I_0)|\le 3$ and $|G(I_1)|\le 4$. Now by Theorem \ref{3GEN}, $\sdepth(I_0)\ge (n-1)-1$, and by induction hypothesis $\sdepth(I_1)\ge (n-1)-2$. Therefore, by  \cite[Proposition 3.3]{arxiv.0712.2308}, $\sdepth(I)\ge \min\Set{\sdepth(I_0),\sdepth(I_1)+1}\ge n-2$.
\end{proof}

To conclude, we ask the following question for squarefree monomial ideals.

\begin{ques}
    \label{question1}
    Let $I$ be an $m$-generated squarefree monomial ideal in $S=K[x_1,\dots,x_n]$. Is it true that $\sdepth(I)\ge n-\floor{\frac{m}{2}}$?
\end{ques}

\section*{Acknowledgements}
    The author would like to thank Bernd Ulrich for his help and encouragement during the preparation of this paper. He also thanks the referee for the helpful suggestions.

\bibliographystyle{plain}

\end{document}